\documentclass[12pt]{article}
\usepackage{amsmath}
\usepackage{amssymb}
\usepackage{amscd}
\usepackage{amsfonts}
\usepackage[all]{xy}
\usepackage[a4paper]{geometry}
\usepackage{mathrsfs}
\usepackage{amsthm}

\newcounter{theorem}
\newcounter{corcountrer}

\newcounter{problemcounter}
\newcounter{remarkcounter}

\newtheorem{corollary}[corcountrer]{Corollary}

\newtheorem*{theorem_without_counter}{Theorem}

\newtheorem{prob}[problemcounter]{Problem}
\newtheorem{remark}[remarkcounter]{Remark}

\newcommand{\refm [1]} {(\ref{#1})}
\newcommand{\PP}{\mathbb{P}}

\newcommand{\ber}{\begin{eqnarray}}
\newcommand{\ena}{\end{eqnarray}}
\newcommand{\rules}{\ \mbox{\rule{.5mm}{4mm}}\ }

\newcommand{\be}{\begin{equation}}
\def\la{\label}

\newcommand{\nin}{\noindent}
\newcommand{\non}{\nonumber}

\def\qed{\hfill \vrule height1.3ex width1.2ex depth-0.1ex}
\def\bbb#1{{\rm I\mkern-3.5mu #1}} \def\bba#1#2{{\rm #1\mkern-#2mu\fudge #1}}

\newcommand{\en}{\end{equation}}

\newcommand{\om}{\Omega_N}
\def\reel{\hbox{{\rm R}\kern-1em\hbox{{\rm I} }}}
\def\relatif{\ \hbox{{\rm Z}\kern-.4em\hbox{\rm Z}}}
\def\nat{\hbox{{\rm N}\kern-1em\hbox{{\rm I} } }}
\def\comp{\hbox{{\rm C}\kern-.55em\hbox{{\rm I} } }}
\def\smallcomp{\hbox{\fiverm C}\kern-.35em{\hbox{\fiverm I}}}

\def\fudge{\mathchoice{}{}{\mkern.5mu}{\mkern.8mu}}
\def\bbc#1#2{{\rm \mkern#2mu\vbar\mkern-#2mu#1}}
\def\bbb#1{{\rm I\mkern-3.5mu #1}} \def\bba#1#2{{\rm #1\mkern-#2mu\fudge
#1}}
\def\bb#1{{\count4=`#1 \advance\count4by-64 \ifcase\count4\or\bba
A{11.5}\or \bbb B\or\bbc C{5}\or\bbb D\or\bbb E\or\bbb F \or\bbc
G{5}\or\bbb H\or \bbb I\or\bbc J{3}\or\bbb K\or\bbb L \or\bbb
M\or\bbb N\or\bbc O{5} \or \bbb P\or\bbc Q{5}\orrrr \bb R\or\bbc
S{4.2}\or\bba T{10.5}\or\bbc U{5}\or    \bba V{12}\or\bba
W{16.5}\or\bba X{11}\or\bba Y{11.7}\or\bba Z{7.5}\fi}}
\def\rat{\hbox{{\rm Q}\kern-.70em\hbox{{\rm I} } }}

\def \P {\bbb P}

\title{ Coagulation processes with Gibbsian time evolution}

\author{{\bf Boris L. Granovsky}
\thanks{E-mail: mar18aa@techunix.technion.ac.il} \\
Department of Mathematics, Technion-Israel Institute of Technology,\\
Haifa, 32000, Israel \\ \quad {\bf Alexander V. Kryvoshaev}
\thanks{E-mail: alexkriv@tx.technion.ac.il}\\
Department of Mathematics, Technion-Israel Institute of Technology,\\
Haifa, 32000, Israel.}
\begin{document}
\maketitle \vskip 5cm \nin American Mathematical Society 2000
subject classifications.

\nin Primary-  82C23; Secondary-60J27, 05A18.

\nin Keywords and phrases:   Stochastic processes of coagulation,
Time dynamics, Gibbs distributions.
 \newpage
 \pagestyle{myheadings}

\begin{center}
{\bf Abstract}

\end{center}
We prove that a stochastic process of pure coagulation has at any time $t\ge 0$  a time
dependent Gibbs distribution if and only if the rates $\psi(i, j)$ of single coagulations are of
the form $\psi(i; j) = if (j) + jf(i),$ where $f$ is an arbitrary nonnegative function on the
set of positive integers. We also obtain a recurrence relation for weights of these Gibbs
distributions that allow to derive the general form of the solution and the explicit solutions
in three particular cases of the function $f$. For the three corresponding models, we study
the probability of coagulation into one giant cluster, by time $t>0.$

\newpage
\section{Process of pure coagulation:\ formulation\\ of the model. Objective and the context\\ of the paper.}

We consider a standard model of stochastic coagulation (see e.g.\cite{DGG}), viewed as
a time continuous Markov chain on the finite set $\Omega_N$ of
partitions $\eta$ of a given  integer $N:$
$$\om=\{\eta=(n_1,\ldots,n_N):\sum_{k=1}^N kn_k=N\}.$$
In the context considered, $N$ is the total number (=total mass)
of identical particles (molecules, planets, animals etc)
partitioned into clusters  of different sizes, so that
$n_k$ is the number of clusters of size $k$ in a partition
$\eta\in\om$.
 Infinitesimal in time transitions  are  coagulations
   of any two clusters of sizes $i$ and $j$ into one cluster of size $i+j$, resulting in the  state transition
$$\eta\rightarrow \eta^{(i,j)}. $$
Here
$\eta^{(i,j)}\in\Omega_N$ codes the state that is obtained
  from a state $\eta\in \om,$ with $n_i>0, n_j>0$, by a   coagulation of any two specific clusters  of sizes $i$ and $j$.
In the sequel such coagulations are called single and their rates are denoted $\psi(i,j)$.
  The following  assumptions on the  rates $\psi(i,j)$  and on the induced  rates
  $K(\eta\rightarrow \eta^{(i,j)})$ of state  transitions,  describe the  class of
coagulation processes ($CP$'s) considered:
\begin{itemize}
\item The function $\psi(i,j)$ is nonnegative,  symmetric
   in $i,j$ and is not dependent on $N$;
\item The rate    $K(\eta\rightarrow \eta^{(i,j)})$ is equal to
the sum of rates of all
  single coagulations $\psi(i,j)$ of  $n_i>0$ groups  of size $i$
   with $n_j>0$ groups of size $j $ each one, so that \ber
  \non K(\eta\rightarrow \eta^{(i,j)})&=& n_in_j\psi(i,j),\quad i\neq j,\quad 2\leq
  i+j\leq N ,\\ \non K(\eta\rightarrow \eta^{(i,i)})&=&\frac{ n_i(n_i-1)}{2}\ \psi(i,i),\quad 2\leq
  2i\leq N.
  \ena
\end{itemize}

 An interpretation of the above expression  for the coagulation kernel $K$ in
terms of the kinetics of droplets of different masses  is given in
\cite{mar}. $CP$'s with rates of  state transitions of
  the above form  are naturally  called mean-field models,
  meaning  that at any state $\eta\in \Omega_N,$ any cluster can coagulate with any other. Hence,  given $N$,
the distribution of a $CP=CP(N)$ at any time $t\ge 0$ is uniquely specified by the initial distribution on the set $\Omega_N$ and the  rates $\psi(i,j).$

The history of $CP$'s goes back to 1918 when Smoluchowski formulated a deterministic
version of the model of pure coagulation of molecules in chemical kinetics in his seminal
paper. Because of the ever growing field of applications and a rich probabilistic context,
the study of a variety of versions of the model continues to be a hot topic in the theory
of stochastic processes. Marcus \cite{mar} was apparently the first to formulate the stochastic
version of a $CP$. A particular case of a pure coagulation, when $\psi(i; j) = a(i+j)+b, a, b \ge 0$
is known in the literature as the Marcus-Lushnikov stochastic process, while in \cite{ald1} as well
as in some other papers, the name is given to all stochastic $CP$'s with rates of single
coagulations of the form $N^{-1}\psi(i; j),$ with an arbitrary $\psi(i,j).$ It is important to point
out that, in contrast to the latter Marcus-Lushnikov process, the basic assumption of our setting is
that the rates of single coagulations  do not depend on $N$.
The equilibria of some reversible models with rates of coagulation and fragmentation
depending on N were studied in \cite{dong1},\cite{dong2}.

Let  $X_N^{(\rho)}(t),\ t\ge 0$ denote a $CP(N)$ starting from an
initial distribution $\rho$ on $ \om.$ The objective of our study
is  the probability distribution (=transition probability) $p(\eta,\rho;t)$ of the process, which is
\be\PP(X_N^{(\rho)}(t)=\eta),\quad \eta\in\om,\quad t\ge 0.\la{xup}
\en
In the sequel we refer in more details to the literature related to the aforementioned
objective.

We describe now the organization of the paper. In Section 2 we formulate our main
result which is the characterization of $CP$'s possessing a probability distribution \refm[xup] of
Gibbsian type at any time $t\ge 0$. As corollaries of our theorem we derive the general
form of weights of the aforementioned Gibbs distributions and the explicit expressions for
the weights for three particular models of $CP$'s. We also analyze the behaviour in time
of some important functionals of the models. In the final section, Section 3, we explain
 the name ``Gibbsian'' given to the class of distributions considered, describe
the linkage to coagulation-fragmentation processes on set partitions and indicate the
relation of Gibbsian  distributions to the theory of random combinatorial structures.

\section{   Main result}
Recall that  $X_N^{(\rho)}(t),\ t\ge 0$ denotes a $CP(N)$ starting from an
initial distribution $\rho$ on $ \om.$  Our goal is  to identify $CP$'s
with a probability distribution \refm[xup] of the form
$$  p(\eta,\rho;t)= C_N(t)
\prod_{k=1}^N\frac{(a_{k,N}(t))^{n_k}}{n_k!},$$ \be
\eta=(n_1,\ldots,n_N)\in \Omega_N, \ t\ge 0.\label{form1}
\end{equation}
Here $C_N(t), \ t\ge 0$ is  time-dependent partition function
and $a_{k,N}(t)\ge 0, \ t\ge 0$ are time-and $N$-dependent
weights. The distributions in the right hand side of \refm[form1] are called
Gibbsian (=canonical Gibbsian). We note that by the above definition the initial distribution $\rho$ is Gibbsian. By virtue of the mass conservation law $\sum_{j=1}^N
jn_j=N,$ tilting transformations of the weights in \refm[form1]
with an arbitrary function $r_N(t)>0, \ t\ge 0:$
$$\tilde{a}_{k,N}(t)=\big(r_N(t)\big)^k a_{k,N}(t),\quad k=1,\ldots,N,\quad t\ge 0$$
and the induced transformation $\tilde{C}_N(t)=(r_N(t))^{-N}C_{N}(t)$ of the partition function, result in different representations  of the generic Gibbs distribution. In view of this fact we assume in
the rest of the paper that the weights $a_{k,N}(t)$ in \refm[form1] are such that
the partition function does not depend on $t\ge
0$, i.e. $C_N(t)=C_N,\ t\ge 0.$ The latter assumption appears to be of great help for our discussion below.\\
 Our main result is the following characterization.
\begin{theorem_without_counter}\label{gibs_thrm}
The probability distributions $p(\eta,\rho;t)$ of  a $CP$ $X_N^{(\rho)}(t), \ t\ge 0$ are  of a
Gibbsian form \refm[form1] if and only if the following three conditions are satisfied:

(i) The  initial distributions $\rho$ on $\om$ are Gibbsian with
arbitrary weights \ $a_{k,N}=a_{k,N}(0)\ge 0,\ k=1,\ldots,N,\ N\ge 1.$

(ii) The rates of single coagulations are of the form  \be \psi(i,j) =
if(j) + jf(i),\quad  i,j: 1\le i+j\le N,\la{forma1}\en where $f$  is an arbitrary nonnegative function on the set of positive integers.

 (iii) The weights  $ a_{k,N}(t)$ are defined recursively by
$$ a_{1,N}(t)=e^{-(N-1)f(1)t}a_{1,N},\quad t\ge 0,$$
$$ a_{k,N}(t)= \int_0^t\frac{\sum_{i+j = k}
a_{i,N}(u)a_{j,N}(u)(if(j)+jf(i))}{2}\ e^{-(N-k)f(k)(t-u)}du +
a_{k,N}e^{-(N-k)f(k)t},$$\be k=2,\ldots,N,\quad t\ge 0,\la{srec}
\en where $a_{k,N}\ge 0$ are constants implied by the
initial distribution $\rho$ in (i).
\end{theorem_without_counter}
\textbf{Proof.}
Our plan is to show that the assertions  $(i)-(iii)$ of  Theorem are implied by
the  Kolmogorov forward equations. Suppose $\eta_{(i,j)}$ is the state that is obtained from a state
$\eta\in \Omega_N,$ with $n_{i+j} > 0$, by a fragmentation of some cluster of size $i + j \ge 2$ into two
clusters of sizes $i$ and $j.$ Then the equations read as follows:
$$
\frac{d}{d t}p(\eta,\rho;t) = -p(\eta,\rho;t)
\Big(
\sum_{1\leq i<j \leq N} n_in_j \cdot \psi(i,j) +
\sum_{i = 1}^N \frac{n_i(n_i -1)}{2}\cdot \psi(i,i)
\Big) +
$$
$$
+ \sum_{1\leq i<j \leq N} p(\eta_{(i,j)},\rho;t)\cdot
(n_i+1)(n_j+1)\cdot \psi(i,j) + $$\be \sum_{i = 1}^N
p(\eta_{(i,i)},\rho;t) \cdot \frac{(n_i+1)(n_i+2)}{2} \cdot \psi(i,i)
.\la{kolm1}
\en
First, we assume that \refm[form1] holds.
 Substituting \refm[form1] with weights $a_{k,N}(t)$, such that  the partition function does not depend on $t\ge 0,$ gives
$$
\frac{d}{d t}p(\eta,\rho;t) = -p(\eta,\rho;t)
\Big(
\sum_{1\leq i<j \leq N} n_in_j \cdot \psi(i,j) +
\sum_{i = 1}^N \frac{n_i(n_i -1)}{2}\cdot \psi(i,i)
\Big) +
$$
$$
+ p(\eta,\rho;t) \Big( \sum_{1\leq i<j \leq N}
\frac{a_{i,N}(t)a_{j,N}(t)}{a_{i+j,N}(t)}\frac{n_{i+j}}{(n_i+1)
(n_j+1)} (n_i+1)(n_j+1)\cdot \psi(i,j) +
$$
$$
+\sum_{i = 1}^N
\frac{a_{i,N}^2(t)}{a_{2i,N}(t)}\frac{n_{2i}}{(n_i+1)(n_i+2)}
\frac{(n_i+1)(n_i+2)}{2} \cdot \psi(i,i) \Big),
$$
which can be written as
$$
\frac{\frac{d}{d t}p(\eta,\rho;t)}{p(\eta,\rho;t)}=
\sum_{1\leq i,j \leq N}
\frac{a_{i,N}(t)a_{j,N}(t)}{2a_{i+j,N}(t)}n_{i+j}\cdot \psi(i,j)
 - \Big(\sum_{1\leq i<j \leq N} n_in_j \cdot \psi(i,j) +
\sum_{i = 1}^N \frac{n_i(n_i -1)}{2}\cdot \psi(i,i) \Big),
$$
\be \eta\in \Omega_N, \quad t\ge 0.\la{durb}\en In view of \refm[form1] and the   time independence of the partition function, we
proceed as
$$\sum_{k=1}^N \frac{a_{i,N}^{\prime}(t)}{a_{i,N}(t)}n_i
=  \sum_{1\leq i,j \leq N}
\frac{a_{i,N}(t)a_{j,N}(t)}{2a_{i+j,N}(t)}n_{i+j}\psi(i,j) -
\Big(\sum_{1\leq i<j \leq N} n_in_j \cdot \psi(i,j) +
 \sum_{i= 1}^N \frac{n_i(n_i -1)}{2}\cdot \psi(i,i) \Big).
$$ $$ \eta\in \Omega_N, \quad t\ge 0.$$
Finally, we have
$$
\sum_{k=1}^N \frac{a_{i,N}^{\prime}(t)}{a_{i,N}(t)}n_i = \Big( \sum_{k
= 2}^{N} \frac{\sum_{i+j = k}
a_{i,N}(t)a_{j,N}(t)\psi(i,j)}{2a_{k,N}(t)}n_k \Big) - \frac{1}{2}
\Big(\sum_{1\leq i,j \leq N} n_in_j \cdot \psi(i,j) -
 \sum_{i= 1}^N \psi(i,i) \cdot n_i \Big).
$$$$ \eta\in \Omega_N, \quad t\ge 0.$$
We now rewrite the last equation as \be  \frac{1}{2}
\Big(\sum_{1\leq i,j \leq N} n_in_j \cdot \psi(i,j)\Big)
=\sum_{k=1}^N A_{k,N}(t)n_k,\ \quad t\ge 0,\ \text{for all }\
\eta=(n_1,\ldots, n_N)\in \Omega_N,\la{hat}\en where
\be A_{k,N}(t)=\frac{\sum_{i+j = k}
a_{i,N}(t)a_{j,N}(t)\psi(i,j)}{2a_{k,N}(t)}-\frac{a_{k,N}^{\prime}(t)}{a_{k,N}(t)}+
\frac{1}{2}\psi(k,k),\quad k=1,\ldots,N, \la{yuc}\en assuming that for $k=1$
the sum  in the right hand side  of the last expression equals zero. Applying the mass conservation law,
the equation \refm[hat] conforms to
$$\frac{1}{2}
\Big(\sum_{1\leq i,j \leq N} n_in_j \cdot \psi(i,j)\Big)
=\frac{1}{2N}\sum_{1\leq i,j \leq
N}n_in_j\big(jA_{i,N}(t)+iA_{j,N}(t)\big),\ $$$$\text{for all }\
\eta=(n_1,\ldots, n_N)\in \Omega_N,\quad  t\ge 0,$$ from which, using the assumed symmetry of the function $\psi(i,j)$ we
 derive the unique form of $\psi:$\be
\psi(i,j)=\frac{1}{N}\big(jA_{i,N}(t)+iA_{j,N}(t)\big), \quad
(n_1,\ldots, n_N)\in \Omega_N,\quad t\ge 0.\la{zor1}\en
Consequently,$$\psi(i,i)=\frac{2}{N}iA_{i,N}(t),\quad t\ge 0,\quad
i\ge 1,$$ which imposes the following necessary and sufficient
conditions on the coefficients $A_{i,N}(t):$ \be
N^{-1}A_{i,N}(t)=:f(i),\quad i\ge 1 \quad \text{does not depend on
$N\ge 1$ and $t\ge 0$}.\la{zor}\en By virtue of this condition and
\refm[zor1], we find \be \psi(i,i)=2if(i),\quad
\psi(i,j)=if(j)+jf(i),\la{psix}\en which proves the necessity of  part $(ii)$  of the
theorem. It is left to demonstrate that the weights $a_{i,N}(t)$
can be found recursively from  \refm[zor] and \refm[yuc]. Firstly,
we recover $a_{1,N}(t):$
$$N^{-1}\big(-\frac{a_{1,N}^{\prime}(t)}{a_{1,N}(t)}+f(1)\big)=f(1),\quad t\ge 0,$$
\be a_{1,N}(t)=a_{1,N}e^{-(N-1)f(1)t},\quad t\ge 0. \la{a1t}\en where
$a_{1,N}(0)=a_{1,N}>0$ is a constant, given by the initial distribution $\rho$.
 For $k\ge 2,$ \refm[zor] takes the form
$$N^{-1}\Big(\frac{\sum_{i+j = k}
a_{i,N}(t)a_{j,N}(t)\psi(i,j)}{2a_{k,N}(t)}-\frac{a_{k,N}^{\prime}(t)}{a_{k,N}(t)}+
kf(k)\Big)=f(k).$$ Observing that the convolution term
$$M_{k,N}(t):=\sum_{i+j = k}
a_{i,N}(t)a_{j,N}(t)\psi(i,j)=\sum_{i+j = k}
a_{i,N}(t)a_{j,N}(t)(if(j)+ jf(i)), \quad t\ge 0$$ does not depend
on $a_{k,N}(t),$ we arrive at the following first order
differential equation for $a_{k,N}(t):$
\be a_{k,N}^{\prime}(t)=\frac{M_{k,N}(t)}{2}-(N-k)f(k)a_{k,N}(t),\quad k= 2,\ldots,N. \en
Solving it we obtain the recurrence relation that conforms to part $(iii)$  of the theorem: \be
a_{k,N}(t)=e^{-(N-k)f(k)t}\Big(
\int_0^t\frac{M_{k,N}(u)}{2}\ e^{(N-k)f(k)u}du + a_{k,N}\Big),\quad
k=2,\ldots,N,\quad t\ge 0.\la{srec1}\en
Now the
validation that, under conditions $(i)-(iii)$, \refm[form1] is the solution of the Kolmogorov equations \refm[kolm1] is obvious.
\qed\\
{\bf Notes.}\ (1) An obvious consequence of the theorem is that  a $CP$ with given rates $\psi(i,j)$ of
form \refm[forma1], starting from a non Gibbsian initial
distribution, does not possess Gibbsian transition probabilities.  In
particular, if the aforementioned process starts from a mixture of
Gibbs distributions, then its transition probabilities \refm[xup] are the  corresponding
mixtures of Gibbs distributions.

(2) The differential recurrence \refm[srec] was originally derived by Lushnikov in \cite{Lu1}, where
the weights $a_{k,N}(t),\ t\ge 0$ are associated with the generating function for the probabilities $p(\eta,\rho;t),\ t\ge 0.$ Buffet and
Pul\'{e} \cite{buf} proved  the core fact that, given $\psi(i,j)=if(j)+jf(i),$ the Kolmogorov equations are solved by Gibbsian distributions $p(\eta,\rho;t)$ with weights satisfying the differential recurrence \refm[srec].
 Our result strengthens the result of \cite{buf} in the following four directions.
$(i)$  We establish  the necessity of the form \refm[forma1] of rates
for Gibbsian transition probabilities \refm[xup]; $(ii)$  The study in \cite{buf} is
limited to Gibbs distributions with partition functions
$C_N=1,\ N\ge 1,$ which substantially restricts the class of
possible initial distributions. For example, under the above restriction, the initial Gibbs
distribution with  constant weights $a_{k,N}=1$ is not allowed,
since in this case  the explicit expression for $C_N$ as a function of $N$ is not
known. So,  the trick with tilting transformation can not be
applied. $(iii)$ Our theorem is proven for the case of a
nonnegative  (rather than positive) function $f$, which enables us to
treat $CP$'s like those of Becker-D\"{o}ring, which are defined in Corollary $(D)$ below;
 $(iv)$ Solving in Corolary $(A)$ below the recurrence relation \refm[srec], we find a general form of the weights
 $a_{k,N}(t).$

(3) In \cite{buf}, p.1047 it was noted that if $N$ is a multiple of an integer $q\ge 1,$ then the
initial distribution $\rho$ concentrated on the partition of $N$ into $\frac{N}{q}$
 groups of size $q$ each one,
is Gibbsian, with
$$a_{q,N} > 0,\quad  a_{k,N} = 0,\quad  N\ge k\neq q. $$
An important particular case of the above measures is given by $q = 1,$ which in chemistry
corresponds to what is called total dissipation. In this case the initial distribution $\rho$ is
concentrated on the partition $(N, 0,\ldots,0)$ and $C_N = \frac{N!}{a_{1,N}^N},$
 where $a_{1,N} > 0$ is arbitrary,
while $a_{k,N} = 0,\  k = 2,\ldots,N.$
A closer look on the above example in \cite{buf} shows that a class of Gibbs initial
distributions encompasses a variety (but not all) of measures concentrated on single partitions
of $N$. We give an example of one of such measures. Let $N = l + m,$ where $0 < l < m$
and $N$ is not divisible by $l$. This secures that $\eta^*=(n_1^*,\ldots, n_N^*)$ with $n_l^* = n_m^* = 1$ is the unique   partition
of $N$ with $n_l>0,\ n_m>0$. Hence, the measure $\rho,$ such that $\rho(\eta^*)=1$ can be viewed as a Gibbs distribution with arbitrary weights $a_{l,N} > 0,\ a_{m,N} > 0,$ all other
weights equal to zero and the partition function $C_N = (a_{l,N}a_{m,N})^{-1}$.

In the rest of the paper we assume that the initial distribution is the total dissipation
with $a_{1,N} = a_{1,N} (0) = 1,\  N \ge  1,$  so that $C_N = N!, \ N \ge  1$.

{\bf Corollaries}

(A) {\bf The general form of weights } $a_{k,N}(t).$

 The recurrence relation
\refm[srec] allows to find a general form of the parameters $a_{k,N}(t)$, for an arbitrary function
$f\ge 0.$ To do this, we need two more notations. Let
$$\eta^{(j)} = (l_1^{(j)},\ldots
l_j^{(j)}): \sum_{s=1}^j s l_s^{(j)}=j$$
be a partition from the set $\Omega_j$ of integer partitions of $j$ and set $$q(\eta^{(j)};N): =
\sum_{s=1}^jl_s^{(j)}(N-s)f(s),$$
for $j = 1,\ldots,N$ and $f\ge 0$ .

{\bf Assertion.} For given $N$ and a function $f\ge 0,$ the solution $a_{k,N}(t)$ of \refm[srec] has the
form
\be a_{k,N}(t) =
\sum_{\eta^{(k)}\in \ \Omega_k} B(\eta^{(k)};N)\exp\Big(-q(\eta^{(k)};N)t\Big), \quad t\ge 0,\quad  k=1\ldots, N, \la{fkkh}\en
where the coefficients $B(\eta^{(k)};N)$ do not depend on $t\ge 0.$

{\bf Proof}\ is done by induction in $ 1 \le k \le N.$ Recalling our assumption that
$$a_{k;N} := a_{k;N} (0) =
\left\{
  \begin{array}{ll}
    1, & \hbox{if }\  k=1 \\
    0, & \hbox{if}\ k=2,\ldots, N,
  \end{array}
\right.
$$
we obtain from \refm[srec]
$$a_{1,N} (t) = e^{-(N-1)f(1)t}, \quad  t\ge 0,$$
which is of the form \refm[fkkh] with $q(\eta^{(1)};N) = (N -1)f(1),\quad  B(\eta^{(1)};N) = 1,\quad  N\ge   1$.
Consequently, we derive from \refm[srec],
$$ a_{2;N}(t) = \frac{
f(1)}
{2(N - 1)f(1) - (N - 2)f(2)}\Big(
e^{-f(2)(N-2)t} - e^{-2 f(1)(N-1)t}\Big),\quad t\ge 0,$$
which is of the  form \refm[fkkh] induced by the two partitions $(2, 0)$ and $(0, 1)$ of $2$. Let now
\refm[fkkh] hold for $a_{j,N},\  j = 1,\ldots,k-1 $. Substituting \refm[srec] for $a_{j,N} (t),\ j = 1,\ldots, k-1$ into \refm[fkkh], the claim for $a_{k,N} (t)$ follows from the fact that for any given $\eta^{(j)}\in \Omega_j,
\eta^{(k-j)}\in \Omega_{k-j}$,
$$
q(\eta^{(j)};N) + q(\eta^{(k-j)};N) = q(\eta^{(k)};N),$$
where $\eta^{(k)}\in \Omega_k$ is the partition of $k$  obtained by merging the two partitions  $\eta^{(j)},\eta^{(k-j)}.$

(B) {\bf Additive rates of coagulation: $\psi(i,j)=(i+j)v,\
 v>0.$}
In this case the probability distribution $p(\eta,\rho;t)$ can be found explicitly. Namely,

$$
p(\eta,\rho;t)=N!e^{-N(r-1)vt}\big(1-e^{-Nvt}\big)^{N-r}N^{-(N-r)}
\prod_{k=1}^N\frac{k^{(k-1)n_k}}{(k!)^{n_k}n_k!}, \quad t\ge 0,
$$
\be (n_1,\ldots,n_N)\in \Omega_{r,N}, \la{A1} \en where $
\Omega_{r,N}=\{\eta=(n_1,\ldots,n_N)\in
\Omega_N:n_1+\ldots+n_N=r,\quad 1\le r\le N\} $ denotes the set of
all partitions of $N$ with $r$ summands (=clusters).

{\bf Proof}. The case considered conforms to \refm[forma1] with
$f(j)=v>0,\ \ j\ge 1.$ We will seek the solution of the recurrence
relation \refm[srec] in the case considered, in the form \be
a_{k,N}(t)=e^{-(N-k)vt}\big(1-e^{-Nvt}\big)^{k-1}v_{k,N},\quad
k=1,\ldots N,\quad N\ge 1, \la{akn1} \en with constants
$v_{k,N}>0$ that will be   determined from \refm[srec]. For $k=1$
 \refm[akn1] reduces to \refm[a1t] with $v_{1,N}= a_{1,N}=1$.
Next, by the
 induction argument, in $k\ge 1$ \refm[srec] gives
$$a_{k,N}(t)=e^{-(N-k)vt}\Big(\int_0^t\frac{vke^{-(2N-k)vu}(1-e^{-Nvu})^{k-2}e^{(N-k)vu}}{2}du\Big)\sum_{i+j=k}v_{i,N}v_{j,N}=$$$$
e^{-(N-k)vt}\big(1-e^{-Nvt}\big)^{k-1}\frac{N^{-1}k}{2(k-1)}\sum_{i+j=k}v_{i,N}v_{j,N}.$$
In view of \refm[akn1], this leads to the following recurrence
relation for $v_{k,N}$:
$$v_{k,N}=\frac{N^{-1}k}{2(k-1)}\sum_{i+j=k}v_{i,N}v_{j,N},\quad k=2,\ldots,N-1.$$
This is a particular case of the recursion (2.33) in \cite{gran},
whose solution is
$$v_{k,N}=N^{-(k-1)}\frac{k^{k-1}}{k!},\quad k=1,\ldots,N$$
(see (2.39) in \cite{gran}).\qed\\
 From   \refm[A1] we find  the probability $p_{coag,N}(t)$ of coagulation into one giant cluster of size $N$ by time $t>0$:
$$p_{coag,N}(t)=(1-e^{-Nvt})^{N-1}\to 1, \ N\to \infty, \quad t>0.$$
This says that a strong gelation phenomenon holds at any time $t>0$, as $N\to \infty$.

We discuss below the probabilistic meaning of the distribution
\refm[A1]. We rewrite \refm[A1] as
\be p(\eta,\rho;t)=\binom{N-1}{r-1}e^{-N(r-1)vt}\big(1-e^{-Nvt}\big)^{N-r}B_{r,N}^{-1}
\prod_{k=1}^N\frac{k^{(k-1)n_k}}{(k!)^{n_k}n_k!},\quad t\ge 0,\la{xxx}\en
$$(n_1,\ldots,n_n)\in \Omega_{r,N},$$ where $B_{r,N}^{-1}=N(r-1)!(N-r)!N^{-(N-r)}$.
 It was noted in \cite{gran}  that the sequence of weights
$\frac{k^{k-1}}{k!},\ k\ge 1$ is a particular case of weight
sequences satisfying the Gnedin-Pitman  condition of
exchangeability of Gibbs set partitions \cite{gnpit}.   It is easy to see that in \refm[xxx] the
term
$$\binom{N-1}{r-1}e^{-N(r-1)vt}\big(1-e^{-Nvt}\big)^{N-r}$$
expresses the transition probabilities of a pure death process with
rates
$\mu_{r,N}=(r-1)vN,\quad r=1,\ldots,N,$ as in (2.17) of
\cite{gran}. Note that the aforementioned transition probabilities are
binomial distributions with time-dependent probabilities of
success, while   $B_{r,N}$ in \refm[xxx] is the $(N,r)$ partial Bell polynomial (see (2.37) in \cite{gran} for more details).  In the
context of our model $(B_{r,N})^{-1}$  serves as the partition
function for the conditional distribution $\PP(X_N^{(\rho)}(t)=\eta\rules \vert X_N^{(\rho)}(t)\vert=r)$ which is  microcanonical Gibbs distribution on $\Omega_{r,N},$ with
weights $\frac{k^{k-1}}{k!}$ that are independent on $N$ and $t$.
In view of this, we conclude from \refm[xxx] that the distribution of the number of clusters
$\vert X_N^{(\rho)}(t)\vert:= n_1(t)+\ldots+ n_N(t)$
at time $t$ is binomial:
\be \PP(\vert X_N^{(\rho)}(t)\vert=r)=\binom{N-1}{r-1}e^{-N(r-1)vt}\big(1- e^{-Nvt}\big)^{N-r},\quad t\ge 0.\la{juvv}\en
This fact was originally proven by Lushnikov (see (49) in \cite{Lu1}). It is interesting that in the case of the  Marcus-Lushnikov $CP$ with the $N$-dependent additive kernel, the distribution of $\vert X_N^{(\rho)}(t)\vert$ is also binomial, but with a different parameter.  The latter distribution
was derived by Aldous in \cite{ald1}, from the interpretation of the process as the vector of sizes of the continium random tree.

(C) {\bf  Multiplicative rates of coagulation: $\psi(i,j)=2ij. $}
Correspondingly, $f(i)=i,\ i\ge 1,$ so that under the assumed
initial distribution, \refm[srec] becomes
 $$a_{1,N}(t)=e^{-(N-1)t}, \ t\ge 0,\quad  N\ge 1,$$
 $$ a_{k,N}(t)=e^{-(N-k)kt} \int_0^te^{k^2u}\sum_{i+j = k}
\big(ia_{i,N}(u)e^{Niu}\big)\big(ja_{j,N}(u)e^{Nju}\big)du, \
k=2,\ldots,N,\quad t\ge 0.$$ Denoting $$ia_{i,N}(u)e^{Niu}=
b_{i,N}(u)$$ we arrive at the following recurrence relation:
$$b_{1,N}(t)=e^t, \ t\ge 0,\ N\ge 1,$$
\be b_{k,N}(t)=ke^{k^2t}\int_0^t e^{-k^2u}\sum_{i+j = k}
b_{i,N}(u)b_{j,N}(u)du, \ k=2,\ldots,N,\quad t\ge 0.\la{srec12}\en
The important fact is that the function $b_{1,N}(t)$ and,
consequently, by virtue of \refm[srec12],  the functions
$b_{k,N}(t)$ do not depend on $N$:
 $b_{k,N}(t)=b_k(t),\ t\ge 0$.
 Thus,  in the case considered,
\be
p(\eta,\rho;t)=N!e^{-N^2t}\prod_{k=1}^N\frac{(b_{k}(t))^{n_k}}{k^{n_k}n_k!},
\ \eta=(n_1,\ldots,n_N)\in \Omega_N, \ t\ge 0, \la{comb}\en where
$b_{k}(t)$ are defined by \refm[srec12]. The right hand side of  \refm[comb] can be viewed as the Gibbs distribution on $\Omega_N$
with the partition function $N!e^{-N^2t}$ and the weights
$\frac{b_k(t)}{k},\ k\ge 1,\ t\ge 0$ not depending on $N.$  This
allows  to employ the known exponential relation (see for
references e.g. \cite{DGG}) between the generating functions
$H(t;x)$ and $V(t;x)$ for the sequences
$h_k(t):=\frac{e^{k^2t}}{k!}, \ k\ge 0, \ t\ge 0$ and
$v_k(t):=\frac{b_k(t)}{k},\quad k\ge 1, \ t\ge 0$ respectively:
\be H(t;x)=e^{V(t;x)}, \ V(t;x)=\sum_{k\ge 1} v_k(t)x^k,\quad H(t;x)=\sum_{k\ge 0}\frac{e^{k^2t}}{k!}x^k,\quad \ t\ge 0.\la{expo}\en
We note that the radius of convergence of the power series
$H(t;x)$ is zero and that  by virtue of  the exponential relation \refm[expo],
the same is true for the series $V(t;x).$ This says that the two
power
 series should be treated as formal ones. From the exponential
relation \refm[expo] it is easy to derive the following recurrence relation between the
sequences $\{v_k(t)\}$ and   $\{h_k(t)\}$: \be h_0(t)\equiv
1,\quad (n+1) h_{n+1}(t)=\sum_{k=0}^n
(k+1)v_{k+1}(t)h_{n-k}(t),\quad n=0,1,\ldots, \quad t\ge
0.\la{expo1}\en It goes without saying that \refm[expo1] and
\refm[srec] are equivalent. However \refm[expo1] is much more
convenient for the study of the asymptotics of $v_k(t),$ as $k\to
\infty$. Based on the fact that the functions
$h_k(t)=\frac{e^{k^2t}}{k!}$ grow very
 rapidly
with $k$ for any fixed $t>0,$ we will
demonstrate that the solution of \refm[expo1] is given by
\be v_k(t)\sim h_k(t),\quad k\to \infty,\quad t>0.  \la{vhtc}\en
Firstly we see from \refm[expo1] that
\be v_{n+1}(t)\le h_{n+1}(t),\quad t\ge 0,\quad n=0,1,\ldots, \la{nfv}\en
This implies
$$(n+1)v_{n+1}\ge (n+1)h_{n+1}(t)- \sum_{k=0}^{n-1}(k+1)h_{k+1}(t)h_{n-k}(t),\quad t\ge 0,$$
and consequently,
$$\frac{v_{n+1}(t)}{h_{n+1}(t)}\ge 1- \sum_{k=0}^{n-1}\frac{(k+1)h_{k+1}(t)h_{n-k}(t)}{(n+1)h_{n+1}(t)},\quad t\ge 0.$$
This together with \refm[nfv] say that for the proof of  \refm[vhtc] one should validate the limit
$$\sum_{k=0}^{n-1}\frac{(k+1)h_{k+1}(t)h_{n-k}(t)}{(n+1)h_{n+1}(t)}= e^{-2nt}\sum_{k=0}^{n-1}\binom{n}{k}e^{-2k(n-k-1)t}\to 0,\quad n\to \infty,  $$
for any $t>0.$
We have,
$$e^{-2nt}\sum_{k=0}^{[n/2]}\binom{n}{k}e^{-2k(n-k-1)t}= e^{-2nt}\Big(1+ \sum_{k=1}^{[n/2]}\frac{1}{k!}\prod_{j=0}^{k-1}(n-j)e^{-2(n-k-1)t}\Big)\le$$
\be e^{-2nt}\Big(1+ \sum_{k=1}^{[n/2]}\frac{1}{k!}ne^{-2(n-k-1)t}\Big)\to 0,\quad n\to\infty.\la{anjv}\en
In view of the relation $\binom{n}{k}=\binom{n}{n-k}$, the same limit can be proven for the sum
$$e^{-2nt}\sum_{k=[n/2]+1}^{n-1}\binom{n}{k}e^{-2k(n-k-1)t},\quad t>0.$$
This completes the proof of \refm[vhtc]. Equipped with this result we are now in a position to study the time dynamics  of clustering of groups of different sizes, as $N\to \infty$ and $t>0$ is fixed. With an obvious abuse of notation, let the random variable $n_{k,N}(t)$ be the number of groups of size $k$ at time $t>0.$  By the known formulae for the functionals of the Gibbs distribution considered (see e.g. \cite{DGG}), we
obtain for any fixed moment of time $t>0$,
$$ En_{k,N}(t)=v_k(t)\frac{h_{N-k}(t)}{h_N(t)}, \quad k=1,\ldots,N.$$
$$ Var\ n_{k,N}(t)= v_k^2(t)\Big(\frac{h_{N-2k}(t)}{h_N(t)}-\big(\frac{h_{N-k}(t)}{h_N(t)}\big)^2\Big)+v_k(t)\frac{h_{N-k}(t)}{h_N(t)},$$
\be cov(n_{k,N}(t),n_{l,N}(t))= v_k(t)v_l(t)\Big(\frac{h_{N-k-l}(t)}{h_N(t)}-\frac{h_{N-k}(t)h_{N-l}(t)}{h_N^2(t)} \Big),
\quad k\neq l=1,2\ldots,N.\la{funct}\en
Substituting $h_k(t)=\frac{e^{k^2t}}{k!},$  it is easy to find
that the three quantities tend to zero as $N\to \infty$ for any fixed $k,l$ and $t>0$. On the other hand,
applying  \refm[vhtc], we conclude that at any time $\quad t>0,$
$$ En_{N,N}(t)\to 1, \quad N\to \infty.$$
This means  that  $p_{coag,N}(t)\to 1,$ as $N\to\infty$ at any moment $ t>0$, which is equivalent to say that a strong form of gelation occurs
during all time evolution of the process.

Unlike this, the Marcus-Lushnikov process with $N$-dependent multiplicative kernel exhibits gelation only after time $t=1$. This fact was proven in \cite{buf}. Also note, that in \cite{buf1} the aforementioned Marcus-Lushnikov process was represented as a random graph process, which  among other things allowed to obtain a version of \refm[comb] for the case considered. The survey
\cite{ald1} enlightens this interesting connection to random graphs.
Finally, we note that a properly time-space rescaled  Marcus-Lushnikov $CP$ with multiplicative kernel converges to a limit process called standard multiplicative coalescence (\cite{ald1}). This fact facilitates the study of the emergence of the  giant component.
Regarding the formulae \refm[funct], it is in order to observe that they are not valid in the general case of
$N$-dependent weights $a_{k,N}(t)$. This can be seen from the derivation of \refm[funct] (see e.g. \cite{DGG}).

(D) {\bf Becker-D\"{o}ring pure coagulation
process.} This $CP$ is a stochastic version of
Becker-D\"{o}ring kinetic equations proposed in 1935, to model a variety of phenomena in
which only coagulations with monomers (=clusters of size $1$) are allowed (see for references
e.g. \cite{BCP}).
 Formally,$$\psi(i,j)=0,\ \text{if}\ \ \min\{i,j\}>1.$$ Clearly,
the rates of the process are of the form \refm[forma1] iff the
function $f$ has  the form \be f(i)= \left\{
                         \begin{array}{ll}
                           0, & \hbox{if }\ \ i>1, \\
                           v>0, & \hbox{if}\ \ i=1,
                         \end{array}
\right. \en which leads to the rates $$
\psi(i,1)=\psi(1,i)=\left\{
\begin{array}{ll}
    iv, & \hbox{if} \ i>1 \\
    2v, & \hbox{if}\ i=1. \\
\end{array}%
\right.
$$
Correspondingly, \refm[srec] takes the  form
$$ a_{1,N}(t)=e^{-(N-1)vt},\quad t\ge 0,$$
$$ a_{k,N}(t)=(k-1)v \int_0^te^{-(N-1)vu}
a_{k-1,N}(u)du, $$\be k=2,\ldots,N,\quad t\ge 0.\la{srec1} \en
From \refm[fkkh] it is not difficult to derive that   in the case considered $a_{k,N}(t)$
is a polynomial of degree $k$ in $z(t)=e^{-(N-1)vt},\ t\ge 0:$
\be a_{k,N}(t)=\sum_{i=0}^k m_{i,k} \cdot e^{-i(N-1)vt},\quad k=1,\ldots,N, \quad t\ge 0.\la{szx}\en
By \refm[srec1], the coefficients $m_{i,k}=m_{i,k,N}$ are  defined recursively by
\be
 m_{i,k} = -\frac{(k-1)m_{i-1,k-1}}{i(N-1)},\quad
i=1,\ldots,k,\quad k\ge 2, \quad m_{0,k} = -\sum_{i=1}^{k}
m_{i,k},\quad k\ge 2, \la{xs} \en with the initial conditions
\be
m_{0,1} =
0,\quad m_{
1,1}=1,\quad
m_{0,2} =a_{2,N}(\infty)= \frac{1}{2(N-1)}, \quad m_{2,2}=-\frac{1}{2(N-1)}.
\la{init}\en
The  recurrence relation in \refm[xs] is validated by the  induction
argument:
$$
 (k-1)v
\int_0^t  e^{-(N-1)vu} \Big( \sum_{i=0}^{k-1}m_{i,k-1} \cdot
e^{-i(N-1)vu} \Big)  du=$$
$$
(k-1)v \int_0^t \Big( \sum_{i=1}^{k} m_{i-1,k-1} \cdot
e^{-i(N-1)vu} \Big)du = (k-1)v
 \sum_{i=1}^{k}
 \frac{m_{i-1,k-1}}{i(N-1)v}(1- e^{-i(N-1)vt})=a_{k,N}(t),$$$$\quad k=2,\ldots,N,\quad t\ge 0.
$$
 Solving \refm[xs], we find  the explicit expressions for the coefficients
$m_{i,k}:$
$$
m_{k-1,k} = 0, \quad k\ge 1,\quad m_{k,k}=\frac{(-1)^{k-1}}{k(N-1)^{k-1}},\quad k\ge 1,
$$
\be m_{i,k}=\frac{(-1)^i(k-1)!}{i!(N-1)^{k-1}(k-i)(k-i-2)!},\quad k-i\ge 2.
\la{coiff}
\en
Consequently, we see that  the
probability of appearance of one giant cluster
of size $N$  decays "almost" exponentially to zero, as $N\to \infty,$ at any time $t>0:$
\be p_{coag,N}(t)=N! a_{N,N}(t)\sim N!m_{0,N}=\frac{(N-1)!}{(N-1)^{N-2}}\sim \sqrt{2\pi} (N-1)^{\frac{3}{2}}e^{-(N-1)},\la{clust}\en
 where the first $\sim$ is due to the last relation in \refm[xs].

In
contrast to the models in (B) and (C), the  Becker-D\"{o}ring $CP$ has a
nontrivial equilibrium distribution (=measure) $\mu_N (\eta)$, which is given by the weights  $a_{k,N}(\infty).$ We have $$
a_{1,N}(\infty)=0,\quad  a_{k,N}(\infty)=m_{0,k}=\frac{k-1}{k(N-1)^{k-1}},\quad k\ge 2,$$
so that
\be
\mu_N(\eta)=N!(N-1)^{-(N-\vert\eta\vert)}\prod_{k=1}^N\Big(\frac{k-1}{k}\Big)^{n_k}\frac{1}{n_k!},
\quad \eta\in \Omega_N,\la{prob}\en where we agree that $0^0=1$ and $\vert\eta\vert=n_1+\cdots+n_N$
is the number of clusters in $\eta$. It follows from  \refm[prob] that
$\mu_N(\eta)=0,$ for all $\eta=(n_1,\ldots,
n_N)\in \Omega_N$  with $n_1>0$.

 Finally, we derive from \refm[prob]  the conditional distribution, given $\vert \eta\vert$,
  at the equilibrium of the process (=microcanonical distribution at equilibrium):
 $$ \mu_N(\eta\rules \vert \eta\vert=l)=\big(B_{l,N})^{-1} \prod_{k=1}^N\Big(\frac{k-1}{k}\Big)^{n_k}\frac{1}{n_k!},$$
 \be \eta=(n_1,\ldots,n_N)\in  \Omega_N:\vert \eta\vert=l\le N, \la{nbxh} \en
where $B_{l,N}=\PP(\vert \eta\vert=l)$ is  the
$(l,N)$ partial Bell polynomial on the set of all partitions of $N$ with $n_1=0,$ induced by the weights $v_k:=
\frac{k-1}{k}$ not depending on $N.$

\section{$CP$'s with canonical and microcanonical Gibbs distributions}

The distribution \refm[form1] can be written as
\be p(\eta,\rho;t)=C_N(t)\exp\Big(-\sum_{i=1}^N(-n_i\log a_{i,N}(t)+ \log(n_i!))\Big),\quad \eta\in \Omega_N,\quad t\ge 0. \la{jrsu}\en
This shows that in the context of statistical physics, \refm[form1] conforms to the canonical Gibbs distribution with time dependent potentials $H_k$ of the following special form:
$$H_k(\eta;t)=0,\quad k\ge 2 \quad \quad H_1(\eta;t)=\sum_{i=1}^N\big( -n_i\log a_{i,N}(t) + \log(n_i!)\big),$$$$ t\ge 0,\quad \eta\in\Omega_N.$$
In the context of stochastic models of coagulation and fragmentation, Gibbs distributions \refm[form1] with weights not depending on $t$ and $N$
emerged as early as  in the 1970's in the works of Kelly (see \cite{Kel}) and Whittle (see \cite{W}), devoted to reversible models of clustering  at equilibrium. Vershik (see e.g. \cite{V1}) who intensively studied such Gibbs distributions in the context of equilibrium models of ideal gas, called them multiplicative measures.
Pitman \cite{Pi} introduced Gibbs processes of pure coagulation/fragmentatation on the state space                       $\Omega_{[N]}, $ which is the set of all partitions of the set $[1,\ldots,N],$ while developing Kingman's theory of exchangeable partitons. Unlike the present paper, Pitman's definition of Gibbs process requires that microcanonical
(rather than canonical) distributions are Gibbsian at any time $t$. Subsequently such processes were extensively studied by
Berestecky and Pitman in \cite{BerP}, the main result of which is the characterization of weights of Gibbs fragmentation processes on $\Omega_{[N]}$. It turned out that the time-reversal of these
processes are $CP$'s with $\psi(i, j) = a(i + j) + b.$ (Note that for $b\neq 0$ the latter rates are not of the form \refm[form1]). In \cite{Gold} Goldschmidt, Martin and Span\`{o} constructed a Gibbs fragmentation process with weights that does not obey the characterization condition in
\cite{BerP}. This became possible because the constructed process does not possess the mean-field
property.
The interplay between the set up, when the state space of a $CP$ is set partitions, and
the set up in the present paper is discussed in more details in \cite{Pi},\cite{BerP},\cite{gran}. In \cite{gran}, as
a development of the idea of \cite{HendSES}, a characterization of coagulation-fragmentation processes, such that the induced birth and death processes $\vert X_N^{(\rho)}(t)\vert$\ are time homogeneous, was established.
 Based on this, a
characterization of coagulation-fragmentation models, possessing time-independent microcanonical Gibbs distributions $\PP( X_N^{(\rho)}(t)= \eta\rules \vert X_N^{(\rho)}(t)\vert=k) $ was obtained. By [11]
and our theorem, the interrelation between $CP$'s with Gibbsian canonical distributions
and the ones with time-independent Gibssian microcanonical distributions is as follows.
(We note that Gibbsian canonical distribution induces a Gibbsian microcanonical distribution, the latter being in general time-dependent).
 The $CP$'s with $\psi(i,j = a(i + j),\  a > 0$ are the only ones that have Gibbsian  canonical
and time-independent Gibbsian microcanical distributions. On the other hand, the $CP$'s with
$\psi(i, j) = a(i + j) + b,\  a \ge 0,\  b > 0$ are the only ones that have time-independent Gibbsian microcanonical distributions and non Gibbsian canonical distributions, while the $CP$'s with
$\psi(i,j)=if(j)+jf(i)$, when $f\neq const$ are the only ones having Gibbsian canoninal
and time-dependent Gibbsian microcanonical distributions.
In the conclusion, we mention a representation of Gibbs distributions \refm[form1] arising
in the field of random combinatorial structures. In the case when the weights $a_{k,N}(t)$
do not depend on $t$ and $N$, the distributions \refm[form1] depict the distributions of vectors
$\eta=(n_1,\ldots,n_N)\in \Omega_N$ of component counts of random combinatorial structures (see
e.g.\cite{ABT},\cite{frgr2},\cite{BerP},\cite{Gran},\cite{Kol}). In this set up $n_k$ stands for the number of non decomposable
components (e.g. cycles in a random permutation) of size $k$. A cornerstone fact in
the theory of random structures is the representation of the aforementioned measures
via the so called conditional relation, which proved to be very useful for problems of
asymptotic enumeration. The version of the conditional relation for time dependent Gibbs
distributions \refm[form1] reads as
$$p(\eta,\rho; t) = \PP\Big(Z_{1,N}(t)=n_1,\ldots,Z_{N,N}(t)=n_N\rules\sum_{k=1}^NkZ_{k,N}(t)=N\big)$$
$$
\eta = (n_1,\ldots n_N)\in\Omega_N,\quad t\ge 0,$$
where $\{ Z_{k,N}(t), k = 1,\ldots,N,\ t\ge 0\}$ is the triangular array of Poisson random variables
with parameters $a_{k,N}(t)$, such that $Z_{k,N}(t),\ k=1,2,\ldots,N $  are independent at any given time
$t > 0$. We note that in our setting,
$$\PP\Big(
\sum_{k=1}^N kZ_{k,N} (t) = N\Big)=(N!)^{-1}exp\big(-\sum_{k=1}^Na_{k,N}(t)\big),\quad t>0.$$

{\bf Acknowledgement} We appreciate constructive suggestions and critical remarks of a
referee.

\end{document}